\documentclass[12pt, article]{amsart}

\usepackage[margin=1.15in]{geometry}
\usepackage{amsmath,amssymb,amsthm}
\usepackage{graphicx}
\usepackage[dvipsnames]{xcolor}
\usepackage{enumitem}
\usepackage[colorlinks=true,linkcolor=MidnightBlue,citecolor=MidnightBlue,urlcolor=MidnightBlue]{hyperref}

\theoremstyle{plain}
\newtheorem{thm}{Theorem}[section]
\newtheorem{lem}[thm]{Lemma}
\newtheorem{prop}[thm]{Proposition}
\newtheorem{cor}[thm]{Corollary}

\newtheorem{mainthm}{Theorem}

\theoremstyle{definition}

\newtheorem{ex}[thm]{Example}
\theoremstyle{remark}
\newtheorem{rem}[thm]{Remark}
\numberwithin{equation}{section}
\newcommand{\Z}{\mathbb{Z}}
\newcommand{\C}{\mathbb{C}}
\newcommand{\R}{\mathbb{R}}
\newcommand{\Q}{\mathbb{Q}}
\DeclareMathOperator{\Res}{Res}
\DeclareMathOperator{\Disc}{Disc}

\renewcommand{\Re}{\operatorname{Re}}
\renewcommand{\Im}{\operatorname{Im}}
\newcommand{\Cas}{\mathcal{C}}
\newcommand{\K}{\mathcal{K}}

\begin{document}

\title[Cassini ovals and translated polynomials]{Cassini ovals, translated polynomials,\\ and the arithmetic of their discriminants and resultants}

\author{Michael O. Oyengo}
% TODO(author): confirm departmental name and postal address.
\address{Department of Mathematics, Maseno University, Private Bag, Maseno, Kenya}
\email{obiero@maseno.ac.ke}

\subjclass[2020]{Primary 30C15, 11R09; Secondary 11B39, 12D10, 26C10, 31A15}
\keywords{Zeros of polynomials, Cassini ovals, lemniscates, Chebyshev polynomials, discriminants, resultants, cyclotomic polynomials, primitive divisors}

\begin{abstract}
For a polynomial $f$, $a>0$ and $b,\omega\in\mathbb{C}$, consider the translate $J^{(a,b,\omega)}f(z)=f(z/a)-\omega f(b/a)$. We show that for $f(x)=(x^{2}-1)^{n}$ and $\omega=\pm1$ every zero of $J^{(a,b,\omega)}f$ lies on the Cassini oval $|z^{2}-a^{2}|=|b^{2}-a^{2}|$. We call these translates \emph{Cassini polynomials} and give their zeros explicitly, prove that the zeros for $\omega=1$ and $\omega=-1$ strictly interlace on each component of the oval, and evaluate their discriminants and resultants in closed form. We then prove  the reduced monic polynomials $Q_{n}=a^{2n}J^{(a,b,\omega)}f/(z^{2}-b^{2})$ factor into homogenized cyclotomic polynomials, satisfying  $\operatorname{Res}(Q_{m},Q_{n})=(b^{2}-a^{2})^{2(m-1)(n-1)}$ for coprime $m,n$, and have discriminants governed by a Lucas sequence. We then prove a comparison principle: if $f$ is even or odd with all zeros in $[-1,1]$, $b\in\mathbb{R}$ and $|\omega|\le1$, then every zero of $J^{(a,ib,\omega)}f$ satisfies $|z^{2}-a^{2}|\le a^{2}+b^{2}$, and the symmetry hypothesis cannot be removed. As an application, every solution of $\det(A(P_{n})-\lambda I)=F_{n+1}$, where $A(P_{n})$ is the adjacency matrix of the path on $n$ vertices and $F_{n+1}$ is a Fibonacci number, satisfies $|\lambda^{2}-4|\le5$. This disproves, for every $n$, a conjecture of Gullerud, Johnson and Mbirika that the real parts of these solutions are unbounded.
\end{abstract}

\maketitle

%%%%%%%%%%%%%%%%%%%%%%%%%%%%%%%%%%%%%%%%%%%%%%%%%%%%%%%%%%%%%%%%%%%%%%%%%
\section{Introduction}\label{sec:intro}
%%%%%%%%%%%%%%%%%%%%%%%%%%%%%%%%%%%%%%%%%%%%%%%%%%%%%%%%%%%%%%%%%%%%%%%%%

\subsection{Zeros of polynomial sequences}
The location of the zeros of a sequence of polynomials $\{P_{n}(z)\}_{n\ge 0}$ is a classical theme joining complex analysis, approximation theory, combinatorics and number theory; see\cite{Marden, RahmanSchmeisser,Szego}. For orthogonal polynomials the zeros are real and interlace in the interval of orthogonality. For the Chebyshev polynomials $T_{n}$ and $U_{n}$ they are explicit and lie in $(-1,1)$~\cite{Mason,Rivlin}. When the sequence of orthogonal polynomials is perturbed, translated or linearly combined, however, the zeros typically leave the real line and accumulate on curves in the complex plane. Describing these curves remains an active area of study.

One of the key results in this area is the theorem of Beraha, Kahane and Weiss~\cite{BKW1,BKW2}. Suppose $P_{n}(z)=\sum_{j=1}^{k}\alpha_{j}(z)\lambda_{j}(z)^{n}$ satisfies a linear recurrence with polynomial coefficients and a mild non-degeneracy condition. Then the limit points of the zeros are the points where two dominant characteristic roots have equal modulus, together with the isolated points where the coefficient of a unique dominant root vanishes. %The resulting \emph{equimodular curves} were studied systematically by Biggs~\cite{Biggs}. They are the source of Sokal's theorem that the chromatic roots of graphs are dense in $\C$~\cite{Sokal}.

Recent developments involve algebraic approaches to seek curves that contain \emph{every} zero of \emph{every} member of the polynomial sequence. Much of the recent progress here is due to Khang Tran and his collaborators. Tran~\cite{Tran2014} proved that if $\sum_{m\ge 0}H_{m}(z)t^{m}=1/D(t,z)$ with $D(t,z)=1+B(z)t+A(z)t^{k}$ and $k\in\{2,3,4\}$, then all zeros of $H_{m}$ lie on an explicit portion of a real algebraic curve. The proof connects root location with the discriminant of $D$ in $t$ and with Ismail's $q$-analogue of the discriminant~\cite{Ismail}. The three-term recurrence $P_{n}+B(z)P_{n-1}+A(z)P_{n-k}=0$ was treated in~\cite{Tran2015}.

Another approach, closest to this current paper, concerns polynomials whose zeros lie on \emph{conics}. Dilcher and Stolarsky~\cite{DS2} found polynomials satisfying non-linear recurrences related to Chebyshev polynomials whose zeros lie close to ellipse-like curves. In his thesis, Oyengo~\cite{OyengoThesis}  showed that for real $\lambda,\kappa\neq 0$ with $\alpha^{2}=\lambda^{2}-\kappa^{2}\neq 0$, the zeros of the \emph{translated Chebyshev polynomials}
\[
T_{n}(z/\alpha)-T_{n}(i\kappa/\alpha)\qquad\text{and}\qquad T_{n}(z/\alpha)-T_{n}(\lambda/\alpha)
\]
lie on the ellipse $x^{2}/\lambda^{2}+y^{2}/\kappa^{2}=1$. They are the images of rotated roots of unity under the linear map $\mathrm{diag}(\lambda,\kappa)$, and in suitable cases they interlace on the ellipse. The analogous translates of $U_{n}$ have their zeros close to, and inside, such ellipses. Independently, Kim and Lee~\cite{KimLee} proved that the zeros of the difference quotient $(T_{n}(z)-T_{n}(a))/(z-a)$ lie on the ellipse $x^{2}/a^{2}+y^{2}/(a^{2}-1)=1$. For the special choice $\lambda=\sqrt5$, $\kappa=1$ the translating values are Lucas and Fibonacci numbers, since $T_{n}(i/2)=\tfrac{i^{n}}{2}L_{n}$.

One might ask the converse problem; given an algebraic curve in the complex plane, can we find a sequence of real polynomials whose zeros lie on this curve? In his thesis, Oyengo showed that for ellipses, one can derive a sequence of translated Chebyshev polynomials whose zeros lie on the ellipses. The same applies for circles which are special ellipses. In this paper, we show that the same is true for Ovals of Cassini, that certain translations of cyclotomic polynomials have their zeros on ovals of Cassini.

\subsection{Cassini ovals}\label{subsec:cassini}
Let $a>0$ and $R>0$. The \emph{Cassini oval} with foci $\pm a$ and product $R$ is the locus of points whose distances to the two foci have product $R$:
\begin{equation}\label{eq:cassini-def}
\Cas(a,R)=\{z\in\C:\ |z-a|\,|z+a|=R\}=\{z\in\C:\ |z^{2}-a^{2}|=R\}.
\end{equation}
In Cartesian, polar and complex coordinates it reads
\begin{equation}\label{eq:cassini-eqs}
\bigl((x-a)^{2}+y^{2}\bigr)\bigl((x+a)^{2}+y^{2}\bigr)=R^{2},\qquad
r^{4}-2a^{2}r^{2}\cos 2\theta+a^{4}=R^{2},\qquad
|z^{2}-a^{2}|=R,
\end{equation}
and expanding the first gives the quartic $(x^{2}+y^{2})^{2}-2a^{2}(x^{2}-y^{2})=R^{2}-a^{4}$. The curve is symmetric about both axes and meets the real axis at $\pm\sqrt{a^{2}+R}$ and, when $R<a^{2}$, also at $\pm\sqrt{a^{2}-R}$. Its shape is governed by the ratio $R/a^{2}$ (Figure~\ref{fig:shapes}):
\begin{enumerate}[leftmargin=2em]
\item $R\ge 2a^{2}$: a single convex oval;
\item $a^{2}<R<2a^{2}$: a single oval with two ``dimples'' at the imaginary axis (peanut shape);
\item $R=a^{2}$: the lemniscate of Bernoulli, with a node at the origin;
\item $0<R<a^{2}$: two disjoint ovals, one around each focus.
\end{enumerate}
Thus $\Cas(a,R)$ is the level set $|P|=R$ of the quadratic polynomial $P(z)=z^{2}-a^{2}$, i.e.\ a \emph{polynomial lemniscate} of degree two. %The geometry of lemniscates was studied by Erd\H{o}s, Herzog and Piranian~\cite{EHP, Ransford,SaffTotik}.

G.~D.~Cassini introduced these curves in the late seventeenth century as candidate planetary orbits competing with Kepler's ellipses~\cite{HirstLloyd,Karatas}. They have since found many applications in Physics, Mathematics Engineering and even in Medicine. For more on other properties of ovals of Cassini, refer to~\cite{Karatas,Lawrence}.%,Weisstein}.

\subsection{Translated polynomials and the Cassini polynomials}
We ask a general question; Let \(f\) be a real polynomial with real zeros and for \(a,b,\omega\in\mathbb{C}\), \(a\neq0\) define the operator \(J^{(a,b,\omega)}\) as
\[J^{(a,b, \omega)}f=f(z/a)-\omega f(b/a).\]

\emph{When do the zeros of \(J^{(a,b,\omega)}f(z)\) lie on, or are controlled by, a single curve depending only on parameters $(a,b,\omega)$?} For Chebyshev polynomials, the curve is an ellipse when \(\omega=\pm1\). In this paper, we show that for the family $f_{n}(x)=(x^{2}-1)^{n}$, the curve is a Cassini oval when \(\omega=\pm1\). We call the resulting translates \emph{Cassini polynomials}, and in this paper, we study their geometry, arithmetic of their discriminants and their resultants. The Cassini oval then turns out to control the translates of every symmetric real-rooted polynomial with zeros in $[-1,1]$. In particular it controls the Chebyshev, Legendre and Gegenbauer translates, including those arising from path graphs.

Throughout, $a>0$ and $b\in\C$ with $b^{2}\notin\{0,a^{2}\}$, with \(\beta=b^{2}-a^{2}\neq 0,\qquad R=|\beta|\). Define the monic polynomials \(G_{n}(z)\) and \(H_{n}(z)\) as
\begin{equation}\label{eq:notation}
G_{n}(z)=(z^{2}-a^{2})^{n}-\beta^{n},\qquad H_{n}(z)=(z^{2}-a^{2})^{n}+\beta^{n},
\end{equation}
so that the $g_{n}=a^{-2n}G_{n}$ and $h_{n}=a^{-2n}H_{n}$ are explicitly given by
\begin{equation}\label{eq:gh-polynomials}
g_{n}(z)=\Bigl(\frac{z^{2}}{a^{2}}-1\Bigr)^{n}-\Bigl(\frac{b^{2}}{a^{2}}-1\Bigr)^{n},\qquad
h_{n}(z)=\Bigl(\frac{z^{2}}{a^{2}}-1\Bigr)^{n}+\Bigl(\frac{b^{2}}{a^{2}}-1\Bigr)^{n}.
\end{equation}
Note that if \(f_{n}(x)=(x^{2}-1)^{n}\) then $g_{n}(z)=f_{n}(z/a)-f_{n}(b/a)$ and $h_{n}(z)=f_{n}(z/a)+f_{n}(b/a)$. Our main results are as follows.

\begin{mainthm}[Geometry; Theorems~\ref{thm:location}, \ref{thm:interlace}]\label{mainA}
For every $n\ge1$ the zeros of $g_{n}$ and $h_{n}$ lie on $\Cas(a,R)$ and are explicitly given by
\[
z=\pm\bigl(a^{2}+\beta e^{2\pi i k/n}\bigr)^{1/2}\quad\text{and}\quad z=\pm\bigl(a^{2}+\beta e^{\pi i(2k+1)/n}\bigr)^{1/2},\qquad 0\le k\le n-1,
\]
respectively. They are simple unless $(1-b^{2}/a^{2})^{n}=1$ $($resp.\ $=-1)$, in which case $0$ is the only multiple zero and it is double. If $R\neq a^{2}$ the zeros of $g_{n}$ and $h_{n}$ strictly interlace on each component of $\Cas(a,R)$. %The normalized zero-counting measures converge weakly to the equilibrium measure of $\K(a,R)=\{|z^{2}-a^{2}|\le R\}$, and $\capa\K(a,R)=R^{1/2}$.
\end{mainthm}

\begin{mainthm}[Arithmetic; Theorems~\ref{thm:disc-g}, \ref{thm:resQ}, \ref{thm:discQ},\ref{thm:zsig}]\label{mainB}
With $Q_{n}=G_{n}/(z^{2}-b^{2})$ and $E_{n}=\sum_{k=0}^{n-1}a^{2k}(a^{2}-b^{2})^{n-1-k}$ identities for resultants and discriminants are given by
\begin{align*}
\Res_{z}(G_{n},H_{n})&=2^{2n}\beta^{2n^{2}},\\
\Disc_{z}(G_{n})&=(2n)^{2n}\beta^{2n(n-1)}\bigl(a^{2n}-(a^{2}-b^{2})^{n}\bigr), \\
\Disc_{z}(H_{n})&=(2n)^{2n}\beta^{2n(n-1)}\bigl(a^{2n}+(a^{2}-b^{2})^{n}\bigr), \\
\Disc_{z}(Q_{n})&=4^{n-1}n^{2n-4}\beta^{2(n-1)(n-2)}E_{n},\\
\Res_{z}(Q_{m},Q_{n})&=\beta^{2(m-1)(n-1)}\ \text{if }\gcd(m,n)=1,\quad\text{and }0\text{ otherwise}.
\end{align*}
If $a^{2},b^{2}\in\Z$ are coprime with $b^{2}<a^{2}$, then for every $n\ge3$, except $(a^{2},b^{2},n)=(1,-1,6)$, there is a prime $p\equiv1\pmod n$ dividing $\Disc_{z}(Q_{n})$ but no $\Disc_{z}(Q_{m})$ with $2\le m<n$.
\end{mainthm}

\begin{mainthm}[Comparison principle; Theorem~\ref{thm:comparison}, Corollary~\ref{cor:pathgraph}]\label{mainC}
Let $f$ be a non-constant polynomial with all zeros in $[-1,1]$ and $f(-x)=\pm f(x)$, let \(a,b\in\R\) with $a>0$ and $b\in\R$. Every zero of \(J^{(a,ib,\omega)}f\) with $|\omega|\le1$ lies in $\K(a,a^{2}+b^{2})$. Every zero of \(J^{(a,b,\omega)}f\), $|b|>a$ and $|\omega|\ge1$ lies outside the interior of $\K(a,b^{2}-a^{2})$. Without the symmetry hypothesis both statements fail. 
\end{mainthm}

Theorem~\ref{mainA} describes the geometry of the zeros of the \(g_{n}\) and \(h_{n}\) polynomials as lying and dense on the ovals of Cassini. The Cassini oval is precisely curve predicted by Tran's discriminant criterion; see Remarks~\ref{rem:BKW} and~\ref{rem:Tran}. What is special here is that the whole zero set, not just its limit, lies on the curve. Theorem~\ref{mainB} is the arithmetic counterpart for the polynomials. The resultant formula is a Cassini analogue of Apostol's cyclotomic resultants and of the Chebyshev resultants of~\cite{JRT,Louboutin}. The discriminant contains a Lucas sequence, which brings in primitive divisors. Theorem~\ref{mainC} proves the Oyengo's earlier conjecture for symmetric polynomials and shows by example that it fails in general.

%\subsection{Organization}
In section~\ref{sec:prelim}, we present some well known results about resultants, discriminants and the Cassini region. In section~\ref{sec:zeros}, we prove results leading to Theorem~\ref{mainA}. In section~\ref{sec:comparison}, we prove the comparison principle and some of its applications. We conclude this paper in section~\ref{sec:disc} with identities for discriminants and resultants involving the Cassini polynomials.% computations. Section~\ref{sec:open} lists open problems.

%%%%%%%%%%%%%%%%%%%%%%%%%%%%%%%%%%%%%%%%%%%%%%%%%%%%%%%%%%%%%%%%%%%%%%%%%
\section{Preliminaries}\label{sec:prelim}
%%%%%%%%%%%%%%%%%%%%%%%%%%%%%%%%%%%%%%%%%%%%%%%%%%%%%%%%%%%%%%%%%%%%%%%%%

\subsection{Resultants and discriminants}
Let $f(x)=a_{m}\prod_{i=1}^{m}(x-\alpha_{i})$ and $g(x)=b_{n}\prod_{j=1}^{n}(x-\beta_{j})$ with $a_{m}b_{n}\neq 0$. We use the following standard facts; see~\cite{GKZ}.
\begin{align}
\Res_{x}(f,g)&=a_{m}^{n}\prod_{i=1}^{m}g(\alpha_{i})=(-1)^{mn}b_{n}^{m}\prod_{j=1}^{n}f(\beta_{j}),\label{eq:res-roots}\\
\Res_{x}(g,f)&=(-1)^{mn}\Res_{x}(f,g),\qquad \Res_{x}(f,g_{1}g_{2})=\Res_{x}(f,g_{1})\Res_{x}(f,g_{2}),\label{eq:res-mult}\\
\Res_{x}(f,g)&=a_{m}^{\,n-\ell}\Res_{x}(f,r)\quad\text{if } g=fq+r,\ \deg r=\ell,\label{eq:res-euclid}\\
\Res_{x}(\lambda f,\mu g)&=\lambda^{n}\mu^{m}\Res_{x}(f,g)\qquad(\lambda,\mu\in\C^{\times}),\label{eq:res-scale}\\
\Disc_{x}(f)&=(-1)^{m(m-1)/2} \frac{1}{a_{m}}\Res_{x}(f,f')=a_{m}^{2m-2}\prod_{i<j}(\alpha_{i}-\alpha_{j})^{2},\label{eq:disc}\\
\Disc_{x}(\lambda f)&=\lambda^{2m-2}\Disc_{x}(f)\qquad(\lambda\in\C^{\times}).\notag
\end{align}
We also use the chain rule of McKay and Wang~\cite{McKayWang}: if $h$ has degree $k$ and leading coefficient $h_{0}$, then
\begin{equation}\label{eq:chainrule}
\Res_{x}\bigl(f(h(x)),g(h(x))\bigr)=\bigl[h_{0}^{\deg f\cdot\deg g}\Res_{x}(f,g)\bigr]^{k},
\end{equation}
unless $h$ is a constant equal to a common root of $f$ and $g$.

The next lemma uses McKay and Wang's chainrule for resultants and computes discriminants of even polynomials. It is the tool we will employ in the calculation of the discriminant formulas in Section~\ref{sec:disc}.

\begin{lem}\label{lem:disc-even}
Let $S\in\C[w]$ be monic of degree $k\ge1$, let $a\in\C$, and put $F(z)=S(z^{2}-a^{2})$. Then
\[
\Disc_{z}(F)=(-4)^{k}\,S(-a^{2})\,\Disc_{w}(S)^{2}.
\]
\end{lem}

\begin{proof}
$F$ is monic of degree $2k$, and $F'(z)=2z\,S'(z^{2}-a^{2})$. Let $\rho_{1},\dots,\rho_{k}$ be the zeros of $S$. The zeros of $F$ are the $2k$ numbers $\pm r_{j}$ with $r_{j}^{2}=a^{2}+\rho_{j}$. By~\eqref{eq:res-roots},
\[
\Res_{z}(F,F')=\prod_{F(r)=0}2r\,S'(r^{2}-a^{2})=2^{2k}\Bigl(\prod_{F(r)=0}r\Bigr)\prod_{j=1}^{k}S'(\rho_{j})^{2}.
\]
Since $F$ is monic of even degree, we have that $F(0)=S(-a^{2})$. Also $\prod_{j}S'(\rho_{j})=\Res_{w}(S,S')=\pm\Disc_{w}(S)$. Finally $\Disc_{z}(F)=(-1)^{2k(2k-1)/2}\Res_{z}(F,F')=(-1)^{k}\Res_{z}(F,F')$.
\end{proof}

\subsection{Cyclotomic polynomials}
For a positive integer \(n\), the \(n-\)th cyclotomic polynomial is the unique irreducible polynomial with integer coefficients that is a divisor of \(x^{n}-1\). Precisely \[\Phi_{n}(x)=\prod _{\stackrel{1\leq k\leq n}{\gcd(k,n)=1}} \left(x-e^{2i\pi\frac{k}{n}}\right).\] Let $\Phi_{e}$ denote the $e$-th cyclotomic polynomial and $\Phi_{e}^{*}(x,y)=y^{\varphi(e)}\Phi_{e}(x/y)$ its homogenization. Then $x^{n}-y^{n}=\prod_{e\mid n}\Phi_{e}^{*}(x,y)$. Put
\[\Psi_{n}(x)=\frac{x^{n}-1}{x-1}=\prod_{\stackrel{ e\mid n}{e>1}}\Phi_{e}(x).
\] We need the classical evaluation
\begin{equation}\label{eq:discPsi}
\Disc_{x}(\Psi_{n})=(-1)^{(n-1)(n-2)/2}\,n^{n-2}\qquad(n\ge2),
\end{equation}
and the following observation. If $\gcd(m,n)=1$, then $\zeta\mapsto\zeta^{n}$ permutes the non-trivial $m$-th roots of unity, so
\begin{equation}\label{eq:Psi-prod}
\prod_{\stackrel{\zeta^{m}=1}{\zeta\neq1}}\Psi_{n}(\zeta)= \prod_{\stackrel{\zeta^{m}=1}{\zeta\neq1}}\frac{\zeta^{n}-1}{\zeta-1}=1\qquad(\gcd(m,n)=1).
\end{equation}

\subsection{The Cassini region}
For $a,R>0$ let $\K(a,R)=\{z:|z^{2}-a^{2}|\le R\}$, the compact region bounded by $\Cas(a,R)$.

\begin{lem}\label{lem:extent}
Let $a,R>0$. For $z\in\K(a,R)$ we have $|z|\le\sqrt{a^{2}+R}$ and
\[
|\Im z|\le\begin{cases}R/(2a),& R\le 2a^{2},\\ \sqrt{R-a^{2}},& R\ge 2a^{2}.\end{cases}
\]
Both bounds are attained on $\Cas(a,R)$.
\end{lem}

\begin{proof}
The first bound follows from $|z|^{2}=|z^{2}|\le|z^{2}-a^{2}|+a^{2}$. Since $\Im z$ is harmonic, its maximum over $\K(a,R)$ is attained on the boundary. Write $z^{2}=a^{2}+Re^{i\theta}$, from $\Im z^{2}=|z^{2}|-\Re z^{2}$ we get $\Im z^{2}=\phi(\cos\theta)$, where $\phi(c)=\sqrt{a^{4}+R^{2}+2a^{2}Rc}-a^{2}-Rc$ for $c\in[-1,1]$. The function $\phi$ is concave, and $\phi'(c)=0$ exactly when $\sqrt{a^{4}+R^{2}+2a^{2}Rc}=a^{2}$, i.e.\ when $c=-R/(2a^{2})$. If $R\le 2a^{2}$ this point lies in $[-1,1]$ and gives $\phi=R^{2}/(2a^{2})$. Otherwise the maximum is at $c=-1$ and equals $2(R-a^{2})$.
\end{proof}

%%%%%%%%%%%%%%%%%%%%%%%%%%%%%%%%%%%%%%%%%%%%%%%%%%%%%%%%%%%%%%%%%%%%%%%%%
\section{Polynomials with zeros on Cassini ovals}\label{sec:zeros}
%%%%%%%%%%%%%%%%%%%%%%%%%%%%%%%%%%%%%%%%%%%%%%%%%%%%%%%%%%%%%%%%%%%%%%%%%
Recall that \(G_{n}(z)=(z^{2}-a^{2})^{n}-\beta^{n},\qquad\text{and}\qquad H_{n}(z)=(z^{2}-a^{2})^{n}+\beta^{n}\) and for $z\in\C$ let $\sqrt{z}$ be the principal square root.

\begin{thm}\label{thm:location}
Let $n\ge1$, $\zeta_{k}=e^{2\pi ik/n}$ and $\eta_{k}=e^{\pi i(2k+1)/n}$ $(0\le k\le n-1)$.
\begin{enumerate}[label=\textup{(\alph*)}]
\item The zeros of $G_{n}(z)$ $($equivalently of $g_{n}(z))$ are $\pm\sqrt{a^{2}+\beta\zeta_{k}}$, and the zeros of $H_{n}(z)$ $($equivalently of $h_{n}(z))$ are $\pm\sqrt{a^{2}+\beta\eta_{k}}$, for $0\le k\le n-1$.
\item All zeros of \(g_{n}(z)\) and \(h_{n}(z)\)  lie on the Cassini oval $\Cas(a,R)$, $R=|b^{2}-a^{2}|$.
\item $G_{n}(z)$ has $2n$ simple zeros unless $(1-b^{2}/a^{2})^{n}=1$. In that case $R=a^{2}$, $\Cas(a,R)$ is a lemniscate of Bernoulli, and $z=0$ $($its node$)$ is a double zero while the others are simple. The same holds for $H_{n}(z)$ with the condition $(1-b^{2}/a^{2})^{n}=-1$.
\end{enumerate}
\end{thm}

\begin{proof}
(a) $G_{n}(z)=0$ if and only if $\bigl((z^{2}-a^{2})/\beta\bigr)^{n}=1$, i.e.\ $z^{2}=a^{2}+\beta\zeta_{k}$ for some $k$. The case of $H_{n}(z)$ is the same with $\zeta_{k}$ replaced by the $n$-th roots $\eta_{k}$ of $-1$.

(b) For every such zero, $|z-a|\,|z+a|=|z^{2}-a^{2}|=|\beta||\zeta_{k}|=R$.

(c) The $n$ values $a^{2}+\beta\zeta_{k}$ are distinct. Each non-zero value has two distinct square roots, which are simple zeros of $G_{n}(z)$. The value $0$ occurs if and only if $\zeta_{k}=-a^{2}/\beta$ for some $k$, i.e.\ if and only if $(-a^{2}/\beta)^{n}=1$, which is equivalent to $(1-b^{2}/a^{2})^{n}=1$. It then occurs for exactly one $k$ and gives the double zero $z=0$. Taking absolute values shows that this forces $R=|\beta|=a^{2}$. The argument for the zeros of  $H_{n}(z)$ is identical.
\end{proof}

\begin{figure}[t]
\centering
\includegraphics[width=\textwidth]{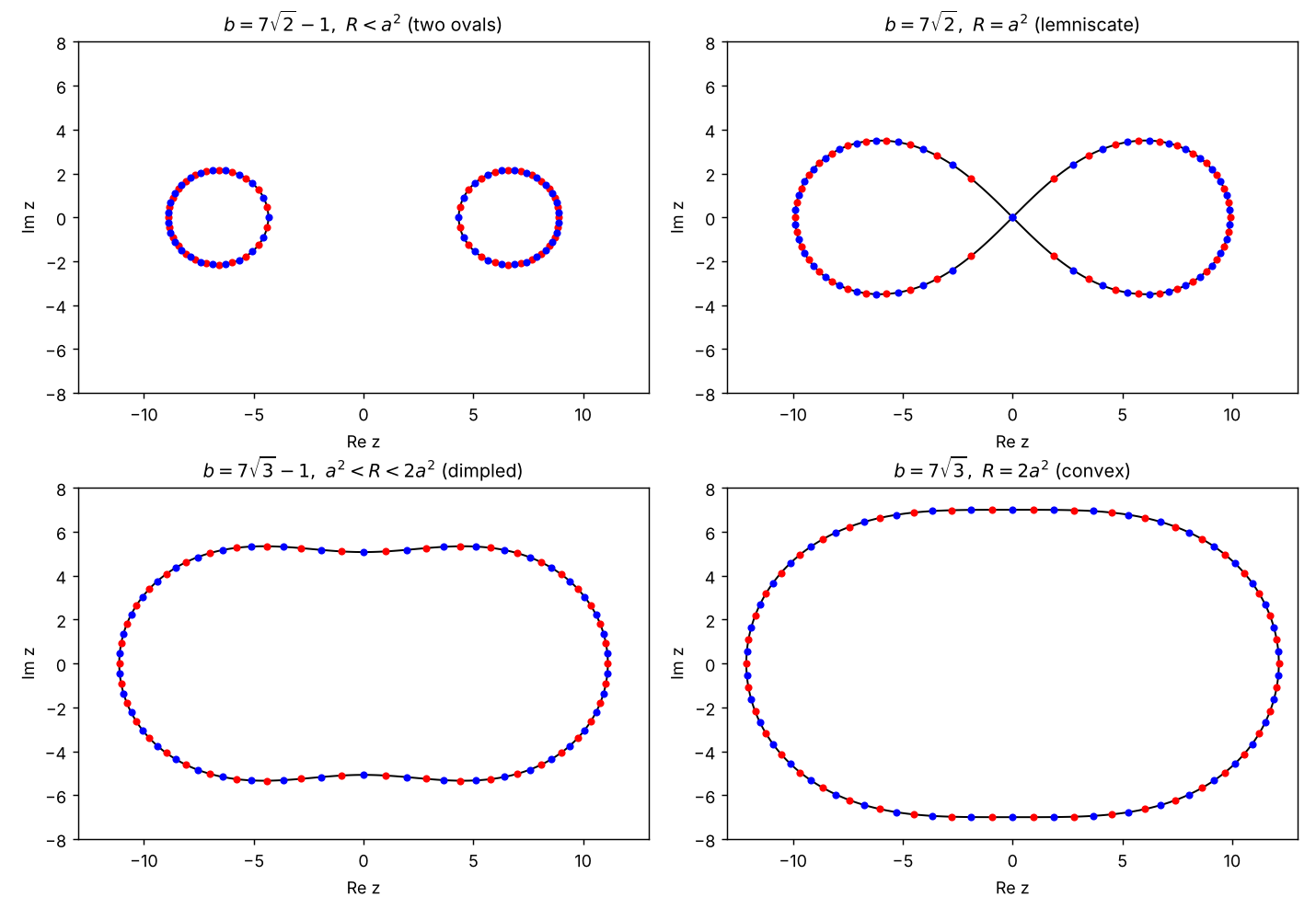}
\caption{Zeros of $g_{12}$ (filled circles) and $h_{12}$ (open squares) for $a=7$ and the indicated values of $b$, on the Cassini ovals $|z^{2}-a^{2}|=|b^{2}-a^{2}|$. Crosses mark the foci $\pm a$. The zeros are equally spaced in the parameter $\arg(z^{2}-a^{2})$ and interlace on each component (Theorem~\ref{thm:interlace}). In the lemniscate case $n=12$ is even and $(1-b^{2}/a^{2})^{12}=1$, so $z=0$ is a double zero of $g_{12}$ (Theorem~\ref{thm:location}(c)).}
\label{fig:shapes}
\end{figure}

\begin{ex}\label{ex:J}
The operator $J^{(a,b,1)}$  applied to $f_{n}(x)=(x^{2}-1)^{n}$, gives $g_{n}(z)$. Hence the zeros of
\(\displaystyle g_{n}(z)=a^{-2n}\bigl((z^{2}-a^{2})^{n}-(b^{2}-a^{2})^{n}\bigr)\)
lie on $\bigl((x-a)^{2}+y^{2}\bigr)\bigl((x+a)^{2}+y^{2}\bigr)=|b^{2}-a^{2}|^{2}$ for every $b\in\C$. The polynomial $J^{(2,i,1)}f_{n}(z)=f_{n}(z/2)-f_{n}(i/2)=4^{-n}\bigl((z^{2}-4)^{n}-(-5)^{n}\bigr)$ has its zeros on $\bigl((x-2)^{2}+y^{2}\bigr)\bigl((x+2)^{2}+y^{2}\bigr)=25$.

For real $a,b>0$, the four regimes of \S\ref{subsec:cassini} become: $b\ge\sqrt3\,a$ (convex); $\sqrt2\,a<b<\sqrt3\,a$ (dimpled); $b=\sqrt2\,a$ (lemniscate); $a<b<\sqrt2\,a$ (two ovals). Figure~\ref{fig:shapes} illustrates all four with $a=7$.
\end{ex}

 \begin{figure}
  \centering
  \includegraphics[width=8cm]{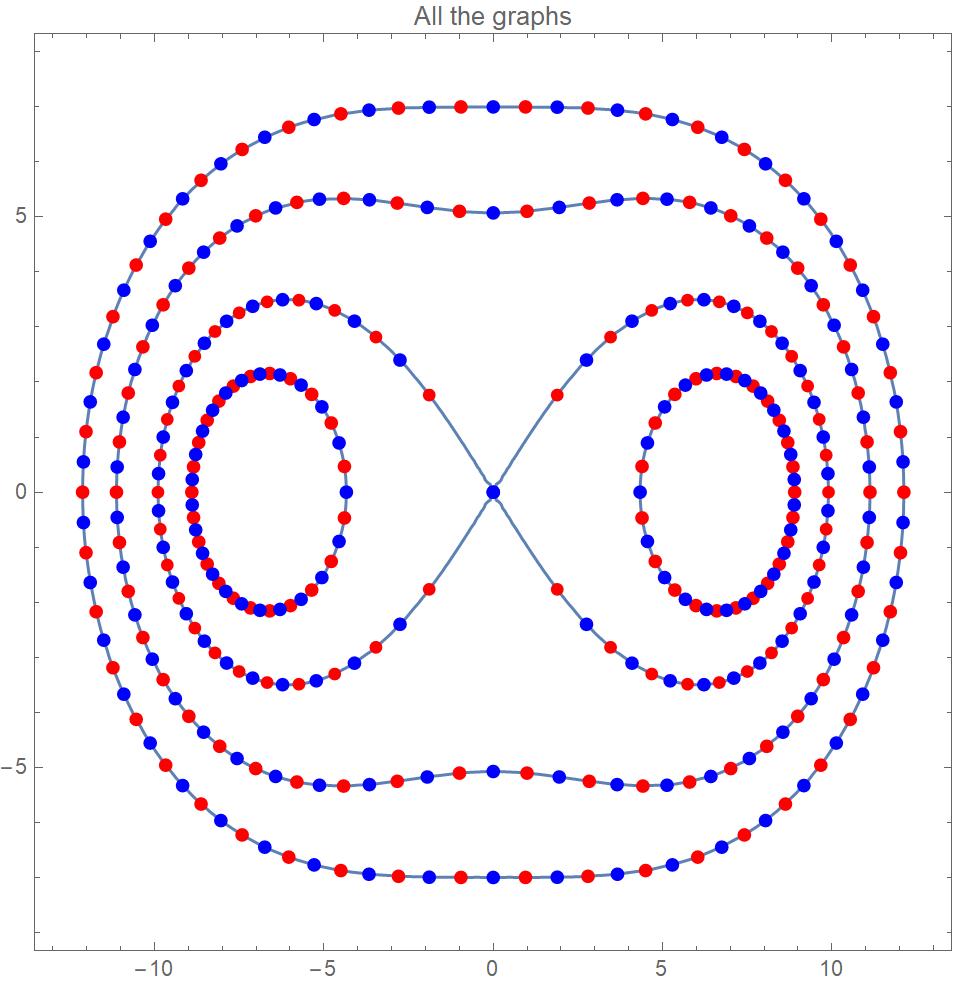}\\
  \caption{All the four curves in one graph.}\label{pic4}
\end{figure}

%\subsection{Parametrization and interlacing}
Recall that \(|z^{2}-a^{2}|=R\) with \(R=|b^{2}-a^{2}|\). When $R\neq a^{2}$ the oval $\Cas(a,R)$ is naturally parametrized by the angle $\theta=\arg(z^{2}-a^{2})$. In this parameter the zeros of $G_{n}$ and $H_{n}$ are equally spaced.

\begin{lem}\label{lem:param}
Let $a,R>0$ with $R\neq a^{2}$ and put $W(\theta)=a^{2}+Re^{i\theta}$.
\begin{enumerate}[label=\textup{(\alph*)}]
\item If $R<a^{2}$, then $\Cas(a,R)$ consists of two Jordan curves $\Gamma_{+}\subset\{\Re z>0\}$ and $\Gamma_{-}=-\Gamma_{+}$. The map $\theta\mapsto\sqrt{W(\theta)}$ is a homeomorphism $\R/2\pi\Z\to\Gamma_{+}$.
\item If $R>a^{2}$, then $\Cas(a,R)$ is a single Jordan curve. There is a continuous $s:\R\to\C$ with $s(\theta)^{2}=W(\theta)$ and $s(\theta+2\pi)=-s(\theta)$, and $s$ induces a homeomorphism $\R/4\pi\Z\to\Cas(a,R)$.
\end{enumerate}
\end{lem}

\begin{proof}
(a) If $R<a^{2}$ then $\Re W(\theta)\ge a^{2}-R>0$, so $\sqrt{W(\theta)}$ is continuous with positive real part. It is injective in $\theta$ modulo $2\pi$, because its square $W(\theta)$ is. Its image is the part of $\Cas(a,R)$ in the right half-plane, and $\Cas(a,R)$ does not meet the imaginary axis because $|(iy)^{2}-a^{2}|=y^{2}+a^{2}>R$.

(b) If $R>a^{2}$, the closed path $W$ winds once around $0$ and never passes through it. Hence a continuous square root $s$ exists along $\R$ and changes sign after one turn. If $s(\theta_{1})=s(\theta_{2})$, then $W(\theta_{1})=W(\theta_{2})$, so $\theta_{2}-\theta_{1}\in2\pi\Z$. If $\theta_{2}-\theta_{1}$ is an odd multiple of $2\pi$, then $s(\theta_{1})=-s(\theta_{1})$, which is impossible since $W\neq0$. Hence $s$ is injective on $\R/4\pi\Z$. Its image is all of $\Cas(a,R)$, since both square roots $\pm s(\theta)=s(\theta),s(\theta+2\pi)$ occur.
\end{proof}

\begin{thm}\label{thm:interlace}
Let $R\neq a^{2}$ and $\theta_{0}=\arg\beta$. In the parameterizations of Lemma~\ref{lem:param}, the zeros of $G_{n}$ correspond to $\theta\equiv\theta_{0}+2\pi k/n$ and the zeros of $H_{n}$ to $\theta\equiv\theta_{0}+\pi(2k+1)/n$. Consequently:
\begin{enumerate}[label=\textup{(\alph*)}]
\item if $R>a^{2}$, the $2n$ zeros of $G_{n}$ and the $2n$ zeros of $H_{n}$ strictly interlace on the single curve $\Cas(a,R)$;
\item if $R<a^{2}$, each of $\Gamma_{\pm}$ carries $n$ zeros of $G_{n}$ and $n$ zeros of $H_{n}$, and they strictly interlace on it.
\end{enumerate}
\end{thm}

\begin{proof}
$W(\theta)=a^{2}+\beta\zeta$ with $|\zeta|=1$ means $\theta\equiv\theta_{0}+\arg\zeta\pmod{2\pi}$. By Theorem~\ref{thm:location}(a) the zeros of $G_{n}$ correspond to $\arg\zeta\in\frac{2\pi}{n}\Z$ and those of $H_{n}$ to $\arg\zeta\in\frac{\pi}{n}+\frac{2\pi}{n}\Z$. These two arithmetic progressions alternate on $\R/2\pi\Z$, and also on its double cover $\R/4\pi\Z$. Homeomorphisms of circles preserve cyclic order, so the claims follow from Lemma~\ref{lem:param}.
\end{proof}

\begin{rem}
When $R=a^{2}$ (the lemniscate) the same argument applies to each lobe. The zeros of $G_{n}$ and $H_{n}$ interlace on each lobe, and one of them may have a double zero at the node $z=0$ (Figure~\ref{fig:shapes}). Theorem~\ref{thm:interlace} also contains an interlacing statement for a single $n$ across different translates. The zeros of $G_{n}$ for $\beta$ and for $\beta e^{i\pi/n}$ interlace, since the second family is $H_{n}$ for the first.
\end{rem}

\begin{rem}\label{rem:BKW}%[Beraha--Kahane--Weiss]
Since $G_{n}=\lambda_{1}^{n}-\lambda_{2}^{n}$ with $\lambda_{1}(z)=z^{2}-a^{2}$ and $\lambda_{2}=\beta$, the sequence $(G_{n})$ satisfies the recurrence $G_{n+1}=(\lambda_{1}+\lambda_{2})G_{n}-\lambda_{1}\lambda_{2}G_{n-1}$ with $G_{0}=0$. Its characteristic roots have equal modulus exactly on $\Cas(a,R)$. The coefficients $\pm1$ never vanish, so Beraha--Kahane--Weiss Theorem ~\cite{BKW1,BKW2} identify $\Cas(a,R)$ as the set of limit points of zeros. Theorem~\ref{thm:location} says much more: every zero of every $G_{n}$ is on the curve.
\end{rem}

\begin{rem}\label{rem:Tran}%[Tran's discriminant criterion]
Let $Q_{n}=G_{n}/(z^{2}-b^{2})$ (see Section~\ref{sec:disc}). Then
\[
\sum_{n\ge1}Q_{n}(z)\,t^{n-1}=\frac1{(1-(z^{2}-a^{2})t)(1-\beta t)}=\frac{1}{1+B(z)t+A(z)t^{2}},
\]
with $B(z)=-(z^{2}+b^{2}-2a^{2})$ and $A(z)=\beta(z^{2}-a^{2})$. This is exactly the setting of Tran's discriminant criterion ~\cite[Theorem~1]{Tran2014}, which places the zeros of $Q_{n}$ on the curve $\Im(B^{2}/A)=0$, $0\le B^{2}/A\le4$. Here $B^{2}/A=2+\lambda+\lambda^{-1}$ with $\lambda=(z^{2}-a^{2})/\beta$. This quantity is real and lies in $[0,4]$ if and only if $|\lambda|=1$, i.e.\ exactly on $\Cas(a,R)$. The $t$-discriminant $B^{2}-4A=(z^{2}-a^{2}-\beta)^{2}$ vanishes only at $z=\pm b$, which is the common zero of all $G_{n}$ removed in $Q_{n}$.
\end{rem}

\begin{rem}%[General lemniscates]
The argument of Theorem~\ref{thm:location} applies verbatim to any non-constant $P\in\C[z]$. The zeros of $P(z)^{n}-P(b)^{n}$ lie on the lemniscate $|P(z)|=|P(b)|$. Cassini ovals are the case $P(z)=z^{2}-a^{2}$. The ellipse results for $T_{n}$ recalled in \S\ref{sec:intro} are of a different nature, because $T_{n}$ is not a power of a fixed polynomial.
\end{rem}

%%%%%%%%%%%%%%%%%%%%%%%%%%%%%%%%%%%%%%%%%%%%%%%%%%%%%%%%%%%%%%%%%%%%%%%%%
\section{A comparison principle and its applications}\label{sec:comparison}
%%%%%%%%%%%%%%%%%%%%%%%%%%%%%%%%%%%%%%%%%%%%%%%%%%%%%%%%%%%%%%%%%%%%%%%%%

An outstanding question is to describe the zeros of \(J^{(a,b,\omega)}f\) for real polynomials \(f\) with zeros in $[-1,1]$. The initial conjecture is that the zeros are either inside or on the Cassini oval \(\Cas(a,R)\). The next theorem proves this under a symmetry hypothesis. Example~\ref{ex:counter} shows that the hypothesis cannot be removed.

\begin{thm}\label{thm:comparison}
Let $f\in\C[x]$ be non-constant with all zeros in $[-1,1]$ and $f(-x)=\pm f(x)$. Let $a>0$.
\begin{enumerate}[label=\textup{(\roman*)}]
\item Let $b\in\R$ and $\omega\in\C$ with $|\omega|\le1$. Every zero $z$ of $f(z/a)-\omega f(i b/a)$ satisfies $|z^{2}-a^{2}|\le a^{2}+b^{2}$. If equality holds, then $z=\pm ib$, or $|\omega|=1$ and $f(x)=C(x^{2}-1)^{m}$.
\item Let $c\in\R$ with $|c|>a$ and $\omega\in\C$ with $|\omega|\ge1$. Every zero $z$ of $f(z/a)-\omega f(c/a)$ satisfies $|z^{2}-a^{2}|\ge c^{2}-a^{2}$.
\end{enumerate}
\end{thm}

\begin{proof}
Because the zero set of $f$ is real, contained in $[-1,1]$ and symmetric about $0$ with multiplicities, we can write
\[
f(x)=C\,x^{\varepsilon}\prod_{j=1}^{m}(x^{2}-t_{j}),\qquad C\neq0,\ \varepsilon\in\{0,1\},\ t_{j}\in[0,1].
\]
Put $u=z^{2}/a^{2}$. Then
\begin{equation}\label{eq:modf}
|f(z/a)|^{2}=|C|^{2}\,|u|^{\varepsilon}\prod_{j=1}^{m}|u-t_{j}|^{2}.
\end{equation}

(i) Let $s=b^{2}/a^{2}\ge0$. Since $(ib/a)^{2}=-s$, formula~\eqref{eq:modf} gives $|f(ib/a)|^{2}=|C|^{2}s^{\varepsilon}\prod_{j}(t_{j}+s)^{2}$. Suppose $|z^{2}-a^{2}|>a^{2}+b^{2}$, i.e.\ $|u-1|>1+s$. For every $t\in[0,1]$,
\begin{equation}\label{eq:tri}
|u-t|\ge|u-1|-(1-t)>(1+s)-(1-t)=t+s\ge0.
\end{equation}
Taking $t=t_{j}$, and $t=0$ for the factor $|u|$, we obtain $|f(z/a)|>|f(ib/a)|\ge|\omega f(ib/a)|$. So $z$ is not a zero. If instead $|u-1|=1+s$, then we do not have a strict inequality in~\eqref{eq:tri}. Equality in the triangle inequality for some $t<1$ forces $u-1\in(-\infty,0]$ with $|u-1|=1+s$, i.e.\ $u=-s$, i.e.\ $z=\pm ib$. Hence, if $z\neq\pm ib$ is a zero on the oval, then $|\omega|=1$, $\varepsilon=0$ and all $t_{j}=1$.

(ii) Let $d=c^{2}/a^{2}>1$. Then $|f(c/a)|^{2}=|C|^{2}d^{\varepsilon}\prod_{j}(d-t_{j})^{2}>0$. Suppose $|z^{2}-a^{2}|<c^{2}-a^{2}$, i.e.\ $|u-1|<d-1$. For every $t\in[0,1]$, $|u-t|\le|u-1|+(1-t)<d-t$. Hence $|f(z/a)|<|f(c/a)|\le|\omega f(c/a)|$, and $z$ is not a zero.
\end{proof}

\begin{rem}
Geometrically, $\{u:|u-t|\le t+s\}$ $(0\le t\le1)$ is a family of discs that increases with $t$ and whose members are all internally tangent at $-s$. The largest one, $t=1$, is the image of $\K(a,a^{2}+b^{2})$ under $z\mapsto z^{2}/a^{2}$. The proof shows that each factor of $f$ is individually dominated. Part~(ii) uses the decreasing family $\{|u-t|\le b-t\}$, whose members are all tangent at $b$.
\end{rem}

\begin{ex}\label{ex:counter}
Symmetry is necessary. Take $f(x)=(x-1)^{2}$, $a=2$, $b=1$. The zeros of $f(z/2)-f(i/2)$ are $z=i$ and $z=4-i$, and $|(4-i)^{2}-4|=|11-8i|=\sqrt{185}>5$. So $4-i$ lies outside $\K(2,5)$.
\end{ex}

\begin{cor}\label{cor:OP}
Let $f$ be any of $T_{n}$, $U_{n}$, the Legendre polynomial $P_{n}$, or a Gegenbauer polynomial $C_{n}^{(\lambda)}$ with $\lambda>-\tfrac12$. For $a>0$ and $b\in\R$, the zeros of $f(z/a)-f(ib/a)$ lie in $\K(a,a^{2}+b^{2})$. For $|c|>a$ the zeros of $f(z/a)-f(c/a)$ lie outside the interior of $\K(a,c^{2}-a^{2})$.
\end{cor}

\begin{proof}
These polynomials have parity $(-1)^{n}$ and all their zeros in $(-1,1)$~\cite{Szego}. Apply Theorem~\ref{thm:comparison} with $\omega=1$.
\end{proof}

For $f=T_{n}$ this is consistent with the ellipse theorem of~\cite{OyengoThesis}. The ellipse $x^{2}/(a^{2}+b^{2})+y^{2}/b^{2}=1$ has foci $\pm a$ and lies in $\K(a,a^{2}+b^{2})$, touching $\Cas(a,a^{2}+b^{2})$ at $\pm ib$.

\begin{cor}\label{cor:pathgraph}
Let $f_{n}(\lambda)=\det(A(P_{n})-\lambda I)$ be the characteristic polynomial of the path graph on $n$ vertices, and let $F_{k}$ denote the Fibonacci numbers. For every $n\ge1$, every complex solution $\lambda$ of $f_{n}(\lambda)=F_{n+1}$ satisfies
\[
|\lambda^{2}-4|\le5,\qquad\text{hence}\qquad |\lambda|\le3\quad\text{and}\quad|\Im\lambda|\le\tfrac54 .
\]
In particular the real parts of the solutions are bounded independently of $n$. This disproves, for every $n$, the conjecture of~\cite{GJM} that these real parts are unbounded (Conjecture~1.1(2) in the numbering of~\cite{DavisFarnsworth}); the case $n\equiv0\pmod4$ was obtained in~\cite{DavisFarnsworth}.
\end{cor}

\begin{proof}
It is well known that $f_{n}(\lambda)=U_{n}(-\lambda/2)$~\cite{GJM}. From $U_{n+1}(x)=2xU_{n}(x)-U_{n-1}(x)$ one checks by induction that $U_{n}(i/2)=i^{n}F_{n+1}$. Hence $F_{n+1}=\omega\,U_{n}(i/2)$ with $\omega=i^{-n}$. The substitution $z=-\lambda$ turns $f_{n}(\lambda)=F_{n+1}$ into $U_{n}(z/2)-\omega\,U_{n}(i/2)=0$. Since $U_{n}$ has parity $(-1)^{n}$ and its zeros lie in $(-1,1)$, Theorem~\ref{thm:comparison}(i) with $f=U_{n}$, $a=2$, $b=1$ and the unimodular constant $\omega$ gives $|z^{2}-4|\le5$, and $\lambda^{2}=z^{2}$. The bounds on $|\lambda|$ and $|\Im\lambda|$ follow from Lemma~\ref{lem:extent} with $a=2$ and $R=5\le2a^{2}$.
\end{proof}

\begin{figure}[t]
\centering
\includegraphics[width=0.92\textwidth]{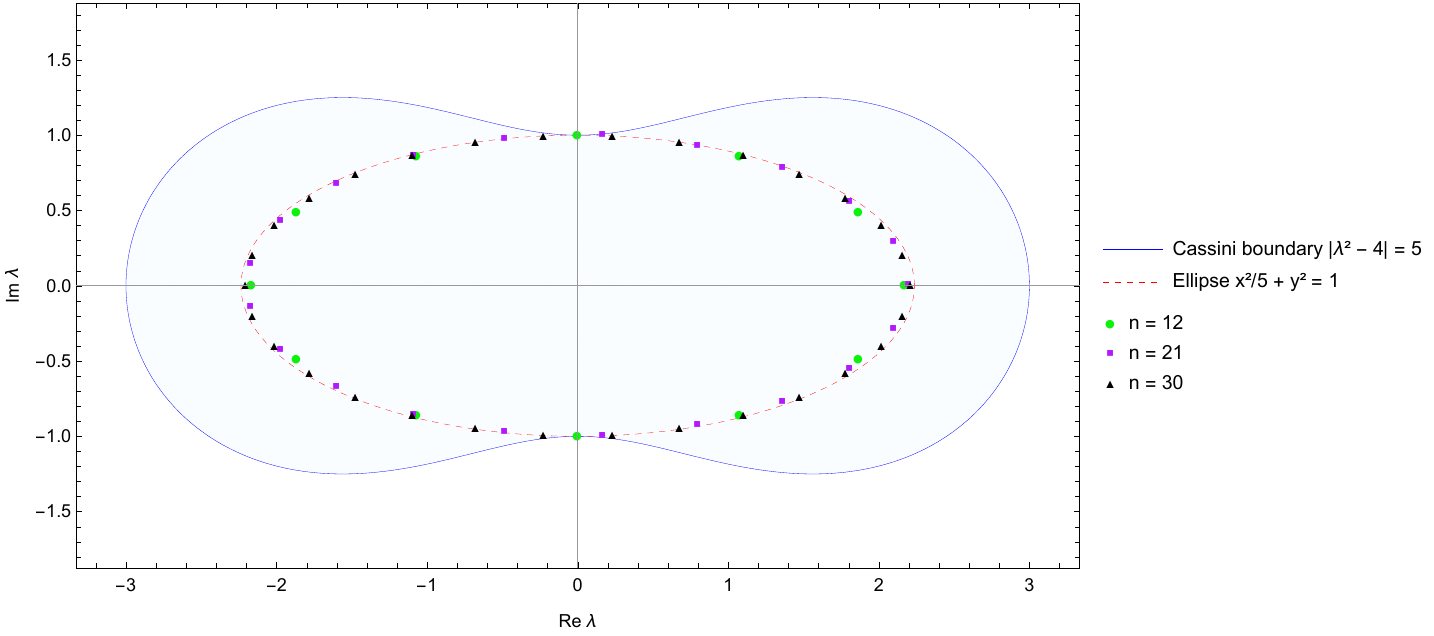}
\caption{Solutions of $f_{P_{n}}(\lambda)=F_{n+1}$ for $n=12,21,30$. They lie in the Cassini region $|\lambda^{2}-4|\le5$ for every $n$ (Corollary~\ref{cor:pathgraph}) and cluster near the ellipse $x^{2}/5+y^{2}=1$ conjectured in~\cite{GJM,DavisFarnsworth}.}
\label{fig:pathgraph}
\end{figure}

\begin{rem}
Numerically the solutions satisfy $|\Im\lambda|\le1$ and $|\Re\lambda|<\sqrt5$ and approach the ellipse $x^{2}/5+y^{2}=1$ (Figure~\ref{fig:pathgraph}), as conjectured in~\cite{GJM,DavisFarnsworth}. Corollary~\ref{cor:pathgraph} gives the uniform but weaker bounds $|\Re\lambda|\le3$ and $|\Im\lambda|\le5/4$. Closing this gap, i.e.\ proving an ellipse inclusion for the translates $U_{n}(z/\alpha)-\omega U_{n}(i\kappa/\alpha)$ with $|\omega|=1$, remains open.%; see Problem~\ref{prob:U}.
\end{rem}

%%%%%%%%%%%%%%%%%%%%%%%%%%%%%%%%%%%%%%%%%%%%%%%%%%%%%%%%%%%%%%%%%%%%%%%%%
\section{Discriminants and resultants}\label{sec:disc}
%%%%%%%%%%%%%%%%%%%%%%%%%%%%%%%%%%%%%%%%%%%%%%%%%%%%%%%%%%%%%%%%%%%%%%%%%

%\subsection{\texorpdfstring{Discriminants of $g_{n}$ and $h_{n}$}{Discriminants of g(n) and h(n)}}

\begin{thm}\label{thm:disc-g}
For $n\ge1$,
\begin{align}
\Disc_{z}(G_{n})&=(2n)^{2n}\,\beta^{2n(n-1)}\bigl(a^{2n}-(a^{2}-b^{2})^{n}\bigr),\label{eq:discG}\\
\Disc_{z}(H_{n})&=(2n)^{2n}\,\beta^{2n(n-1)}\bigl(a^{2n}+(a^{2}-b^{2})^{n}\bigr),\label{eq:discH}
\end{align}
and consequently
\[
\Disc_{z}(g_{n})=a^{-4n(2n-1)}\Disc_{z}(G_{n})\qquad\text{and}\qquad \Disc_{z}(h_{n})=a^{-4n(2n-1)}\Disc_{z}(H_{n}).
\]
\end{thm}

\begin{proof}
Apply Lemma~\ref{lem:disc-even} with $S(w)=w^{n}-\beta^{n}$ and $k=n$. The zeros of $S$ are $\beta\zeta_{k}$, so $S'(\beta\zeta_{k})=n\beta^{n-1}\zeta_{k}^{n-1}$ and
\[
\Disc_{w}(S)^{2}=\Bigl(\prod_{k}S'(\beta\zeta_{k})\Bigr)^{2}=n^{2n}\beta^{2n(n-1)}\Bigl(\prod_{k}\zeta_{k}\Bigr)^{2(n-1)}=n^{2n}\beta^{2n(n-1)},
\]
because $\prod_{k}\zeta_{k}=\pm1$. Moreover $(-4)^{n}S(-a^{2})=4^{n}\bigl(a^{2n}-(-\beta)^{n}\bigr)=4^{n}\bigl(a^{2n}-(a^{2}-b^{2})^{n}\bigr)$. This proves~\eqref{eq:discG}. For~\eqref{eq:discH} take $S(w)=w^{n}+\beta^{n}$; the same computation gives $(-4)^{n}S(-a^{2})=4^{n}\bigl(a^{2n}+(a^{2}-b^{2})^{n}\bigr)$. The formula for $g_{n}=a^{-2n}G_{n}$ and  $h_{n}=a^{-2n}H_{n}$ follow from~\eqref{eq:disc} with $\lambda=a^{-2n}$ and $\deg G_{n}=\deg H_{n}=2n$ respectively.
\end{proof}

The vanishing of~\eqref{eq:discG} is exactly the condition of Theorem~\ref{thm:location}(c). For real $a,b$ it happens only when $b^{2}=2a^{2}$ and $n$ is even, i.e.\ on the lemniscate, when the node $z=0$ is a zero. This fits the geometry.

\begin{prop}\label{prop:resGH}
$\Res_{z}(G_{n},H_{n})=2^{2n}\beta^{2n^{2}}$, and hence $\Res_{z}(g_{n},h_{n})=2^{2n}\beta^{2n^{2}}a^{-8n^{2}}$.
\end{prop}

\begin{proof}
$H_{n}=G_{n}+2\beta^{n}$, so by~\eqref{eq:res-roots} $\Res_{z}(G_{n},H_{n})=\prod_{G_{n}(r)=0}2\beta^{n}=(2\beta^{n})^{2n}$. The normalized version follows from~\eqref{eq:res-scale}.
\end{proof}

%\subsection{\texorpdfstring{The reduced polynomials $Q_{n}$}{The reduced polynomials Q(n)}}
Every $G_{n}$ vanishes at $z=\pm b$, so resultants $\Res(G_{m},G_{n})$ are identically zero. The natural objects are the reduced polynomials
\begin{equation}\label{eq:Qn}
Q_{n}(z)=\frac{G_{n}(z)}{z^{2}-b^{2}}=S_{n}(z^{2}-a^{2}),\qquad S_{n}(w)=\frac{w^{n}-\beta^{n}}{w-\beta}=\beta^{n-1}\Psi_{n}(w/\beta),
\end{equation}
which are monic of degree $2(n-1)$. They are Cassini analogues of the Chebyshev difference quotients of Kim and Lee~\cite{KimLee}. Their generating function was given in Remark~\ref{rem:Tran}, so they satisfy
\[
Q_{n+1}(z)=(z^{2}+b^{2}-2a^{2})\,Q_{n}(z)-\beta(z^{2}-a^{2})\,Q_{n-1}(z),\qquad Q_{0}=0,\ Q_{1}=1.
\]
If $a^{2},b^{2}\in\Z$, then $Q_{n}\in\Z[z]$.

\begin{prop}[Cyclotomic factorization]\label{prop:cyclo}
For $n\ge2$, \(Q_n\) has the cyclotomic factorization;
\[
Q_{n}(z)=\prod_{e\mid n,\ e>1}\Phi_{e}^{*}\bigl(z^{2}-a^{2},\ \beta\bigr),
\]
and the factor with index $e$ has degree $2\varphi(e)$ in $z$. Its zeros are the points of $\Cas(a,R)$ with $\arg(z^{2}-a^{2})-\arg\beta\in\frac{2\pi}{e}(\Z/e\Z)^{\times}$.
\end{prop}

\begin{proof}
Substitute $x=z^{2}-a^{2}$ and $y=\beta$ in $x^{n}-y^{n}=\prod_{e\mid n}\Phi_{e}^{*}(x,y)$ and cancel the factor $\Phi_{1}^{*}=x-y=z^{2}-b^{2}$. The description of the zeros follows from Theorem~\ref{thm:location}(a).
\end{proof}

\begin{rem}\label{rem:quadratic}
The factor with $e=2$ is $\Phi_{2}^{*}(z^{2}-a^{2},\beta)=z^{2}-(2a^{2}-b^{2})$. Its zeros $\pm\sqrt{2a^{2}-b^{2}}$ are points where the oval meets a coordinate axis: the real axis if $b^{2}<2a^{2}$, the imaginary axis if $b^{2}>2a^{2}$. For $a,b\in\Z$ it splits over $\Q$ if and only if $b^{2}+c^{2}=2a^{2}$ for some $c\in\Z$. This conic has infinitely many integer solutions, e.g.\ $(a,b,c,a)=(5,1,7),(13,7,17),(17,7,23),(25,31,17)$. The ratio $R/a^{2}=|b^{2}-a^{2}|/a^{2}$ then records which regime of \S\ref{subsec:cassini} the oval is in. In general, the factorization of $\Phi_{e}^{*}(z^{2}-a^{2},\beta)$ over $\Q$ is a Diophantine question.%; see Problem~\ref{prob:irr}.
\end{rem}

\begin{thm}\label{thm:resQ}
For $m,n\ge2$,
\[
\Res_{z}(Q_{m},Q_{n})=\begin{cases}\beta^{2(m-1)(n-1)}=(b^{2}-a^{2})^{2(m-1)(n-1)},&\gcd(m,n)=1,\\[2pt] 0,&\gcd(m,n)>1.\end{cases}
\]
\end{thm}

\begin{proof}
By~\eqref{eq:Qn} and the chain rule~\eqref{eq:chainrule} with $h(z)=z^{2}-a^{2}$, $h_{0}=1$, $\deg h=2$, we have $\Res_{z}(Q_{m},Q_{n})=\Res_{w}(S_{m},S_{n})^{2}$. The zeros of $S_{m}$ are $\beta\zeta$ with $\zeta^{m}=1$, $\zeta\neq1$. Since $S_{m}$ is monic,
\[
\Res_{w}(S_{m},S_{n})=\prod_{\zeta^{m}=1,\ \zeta\neq1}S_{n}(\beta\zeta)=\beta^{(m-1)(n-1)}\prod_{\zeta^{m}=1,\ \zeta\neq1}\Psi_{n}(\zeta).
\]
If $e=\gcd(m,n)>1$, the factor with $\zeta=e^{2\pi i/e}$ vanishes. If $\gcd(m,n)=1$, the last product equals $1$ by~\eqref{eq:Psi-prod}.
\end{proof}

\begin{rem}
Geometrically, $Q_{m}$ and $Q_{n}$ share a zero if and only if $\gcd(m,n)>1$. In that case, by proposition~\ref{prop:cyclo}, they share exactly the $2\sum_{e\mid\gcd(m,n),\,e>1}\varphi(e)=2(\gcd(m,n)-1)$ zeros of $Q_{\gcd(m,n)}$. When $a^{2},b^{2}\in\Z$ and $\gcd(m,n)=1$, the resultant is the $2(m-1)(n-1)$-th power of the integer $\beta$. This parallels Apostol's theorem that $\Res(\Phi_{m},\Phi_{n})=1$ unless $m/n$ is a prime power~\cite{Apostol}, and the Jacobi-symbol formulas for Chebyshev resultants~\cite{JRT,Louboutin}.%s,Yamagishi}.
\end{rem}

\begin{thm}\label{thm:discQ}
For $n\ge2$,
\[
\Disc_{z}(Q_{n})=4^{n-1}\,n^{2n-4}\,\beta^{2(n-1)(n-2)}\,E_{n},\qquad
E_{n}=\frac{a^{2n}-(a^{2}-b^{2})^{n}}{b^{2}}=\sum_{k=0}^{n-1}a^{2k}(a^{2}-b^{2})^{n-1-k}.
\]
\end{thm}

\begin{proof}
Apply Lemma~\ref{lem:disc-even} with $S=S_{n}$ and $k=n-1$. By~\eqref{eq:Qn} and~\eqref{eq:disc}, $\Disc_{w}(S_{n})=\beta^{(n-1)(n-2)}\Disc_{x}(\Psi_{n})$. By~\eqref{eq:discPsi}, $\Disc_{w}(S_{n})^{2}=\beta^{2(n-1)(n-2)}n^{2(n-2)}$. Finally
\[
(-4)^{n-1}S_{n}(-a^{2})=(-4)^{n-1}\frac{(-a^{2})^{n}-\beta^{n}}{-a^{2}-\beta}=4^{n-1}\frac{a^{2n}-(-\beta)^{n}}{b^{2}}=4^{n-1}E_{n},
\]
since $\beta=b^{2}-a^{2}$.
\end{proof}

For example, $\Disc_{z}(Q_{2})=4(2a^{2}-b^{2})$, the discriminant of the quadratic in Remark~\ref{rem:quadratic}.

\subsection{Primitive prime divisors}
Put $A=a^{2}$ and $B=a^{2}-b^{2}$. Then $E_{n}=(A^{n}-B^{n})/(A-B)$ is the Lucas sequence attached to the pair $(A,B)$. We recall Zsigmondy's theorem~\cite{Zsigmondy}. If $A>B\ge1$ are coprime integers and $n\ge2$, then $A^{n}-B^{n}$ has a prime divisor $p$ that divides no $A^{k}-B^{k}$ with $1\le k<n$ (a \emph{primitive} divisor), except when $n=2$ and $A+B$ is a power of $2$, or $(A,B,n)=(2,1,6)$. A primitive divisor satisfies $p\equiv1\pmod n$, since the order of $A/B$ modulo $p$ is exactly $n$.% For the general theory of primitive divisors of Lucas and Lehmer numbers see Bilu, Hanrot and Voutier~\cite{BHV}.

\begin{thm}\label{thm:zsig}
Let $a^{2},b^{2}\in\Z$ with $b^{2}<a^{2}$, $b^{2}\neq0$ and $\gcd(a^{2},b^{2})=1$. Then for every $n\ge3$ with $(a^{2},b^{2},n)\neq(1,-1,6)$ there is a prime $p\equiv1\pmod n$ such that
\[
p\mid\Disc_{z}(Q_{n})\qquad\text{and}\qquad p\nmid\Disc_{z}(Q_{m})\quad\text{for all }2\le m<n.                                                                                                                       \]
Equivalently, modulo $p$ two zeros of $Q_{n}$ collide, whereas for $m<n$ the zeros of $Q_{m}$ remain distinct modulo $p$.
\end{thm}

\begin{proof}
Put $A=\max(a^{2},a^{2}-b^{2})$ and $B=\min(a^{2},a^{2}-b^{2})$. Both are positive integers, since $b^{2}<a^{2}$ and $a^{2}\ge1$. They are distinct, and $\gcd(A,B)=\gcd(a^{2},b^{2})=1$. Also $A^{n}-B^{n}=\pm\bigl(a^{2n}-(a^{2}-b^{2})^{n}\bigr)$, and $(A,B)=(2,1)$ exactly when $(a^{2},b^{2})=(1,-1)$. By Zsigmondy's theorem there is a primitive prime divisor $p\equiv1\pmod n$ of $A^{n}-B^{n}$. It does not divide $A-B=\pm b^{2}$, so $p\mid E_{n}$, and therefore $p\mid\Disc_{z}(Q_{n})$ by Theorem~\ref{thm:discQ}.

Let $2\le m<n$. We check the factors of $\Disc_{z}(Q_{m})=4^{m-1}m^{2m-4}\beta^{2(m-1)(m-2)}E_{m}$ one by one. Since $p\equiv1\pmod n$ and $n\ge3$, $p$ is odd and $p>n>m$, so $p\nmid 4^{m-1}m^{2m-4}$. If $p\mid\beta=-(a^{2}-b^{2})$, then $p\mid B$ or $p\mid A$, and $p\mid A^{n}-B^{n}$ forces $p\mid\gcd(A,B)=1$, which is impossible. By primitivity $p\nmid A^{m}-B^{m}=\pm b^{2}E_{m}$, so $p\nmid E_{m}$.

The final sentence holds because the discriminant of a polynomial with integer coefficients vanishes modulo $p$ if and only if its reduction modulo $p$ has a repeated factor, or its degree drops. The leading coefficient of $Q_{n}$ is $1$, so the degree does not drop.
\end{proof}

\begin{rem}
The exceptional case $(a^{2},b^{2})=(1,-1)$, i.e.\ $a=1$, $b=i$, has $E_{n}=2^{n}-1$. Here $\Disc_{z}(Q_{6})$ has no new prime, because $2^{6}-1=3^{2}\cdot7$ and $3\mid2^{2}-1$, $7\mid2^{3}-1$. The case $b^{2}>a^{2}$ reduces to primitive divisors of $a^{2n}\pm(b^{2}-a^{2})^{n}$ and can be treated in the same way with the Zsigmondy theorem for $A^{n}+B^{n}$.                                                                                                                                        \end{rem}
%
%%%%%%%%%%%%%%%%%%%%%%%%%%%%%%%%%%%%%%%%%%%%%%%%%%%%%%%%%%%%%%%%%%%%%%%%%
\section{Concluding remarks and open problems}\label{sec:open}
%%%%%%%%%%%%%%%%%%%%%%%%%%%%%%%%%%%%%%%%%%%%%%%%%%%%%%%%%%%%%%%%%%%%%%%%%

In this paper, we have demonstrated that certain translates of real polynomials of the type $f_{n}(x)=(x^{2}-1)^{n}$ have their zeros on Cassini ovals. The results above show that Cassini ovals play for the family $f_{n}(x)=(x^{2}-1)^{n}$ the role that ellipses play for Chebyshev polynomials. Certain questions remain outstanding, and these come from the work of Oyengo~\cite{OyengoThesis, GJM,DavisFarnsworth}. That is: for $\lambda^{2}-\kappa^{2}=\alpha^{2}>0$ and $|\omega|=1$ is there an explicitly defined ellipse on which the zeros of $U_{n}(z/\alpha)-\omega U_{n}(i\kappa/\alpha)$ lie? Oyengo~\cite{OyengoThesis} conjectured that then lie close to the  ellipse $x^{2}/\lambda^{2}+y^{2}/\kappa^{2}\le1$. Corollary~\ref{cor:pathgraph} settles the conjecture $|\Im\lambda|\le1$ of~\cite{GJM} (Conjecture~1.1(3) in~\cite{DavisFarnsworth}) and the bound $|\Re\lambda|\le\sqrt5$ conjectured in~\cite{DavisFarnsworth}.

For symmetric real-rooted $f$ with zeros in $[-1,1]$, Theorem~\ref{thm:comparison} is sharp for $f_{n}(x)=(x^{2}-1)^{n}$. What is the smallest region containing the zeros of $f(z/a)-f(ib/a)$ for all such polynomials $f$ of a fixed degree. What is the analogous region for non-symmetric $f$ (Example~\ref{ex:counter}).

\addcontentsline{toc}{chapter}{\uppercase{References}}
\bibliographystyle{plain} % or unsrt, alpha, etc.
\bibliography{references}

\end{document}